\documentclass{amsart}

\usepackage{amssymb,amsmath}
\usepackage{verbatim,color}
\usepackage{graphicx,xypic}
\usepackage{url}

\newtheorem{theorem}{Theorem}[section]

\newtheorem{proposition}[theorem]{Proposition}
\newtheorem{corollary}[theorem]{Corollary}
\newtheorem{conjecture}[theorem]{Conjecture}

\theoremstyle{definition}
\newtheorem{definition}[theorem]{Definition}
\newtheorem{example}[theorem]{Example}

\theoremstyle{remark}

\numberwithin{equation}{section}

\renewcommand{\AA}{\mathbb{A}}

\newcommand{\CC}{\mathbb{C}}

\newcommand{\FF}{\mathbb{F}}

\newcommand{\PP}{\mathbb{P}}

\newcommand{\RR}{\mathbb{R}}
\newcommand{\ZZ}{\mathbb{Z}}

\newcommand{\into}{\hookrightarrow}     
\def\isomap{{\buildrel \sim\over\longrightarrow}} 
\renewcommand{\iff}{\Leftrightarrow}            

\begin{document}

\title{Hodge Theory in Combinatorics}


\author{Matthew Baker}
\address{School of Mathematics,
          Georgia Institute of Technology, Atlanta GA 30332-0160, USA}
\email{mbaker@math.gatech.edu}
\thanks{The author's research was supported by the National Science Foundation research grant DMS-1529573.  We thank Karim Adiprasito, Eric Katz, and June Huh for numerous helpful discussions and corrections, and Bruce Sagan for additional useful proofreading.}


\date{\today}


\begin{abstract}
If $G$ is a finite graph, a {\em proper coloring} of $G$ is a way to color the vertices of the graph using $n$ colors so that no two vertices connected by an edge have the same color.  (The celebrated four-color theorem asserts that if $G$ is planar then there is at least one proper coloring of $G$ with 4 colors.)  By a classical result of Birkhoff, the number of proper colorings of $G$ with $n$ colors is a polynomial in $n$, called the {\em chromatic polynomial} of $G$.  
Read conjectured in 1968 that for any graph $G$, the sequence of absolute values of
coefficients of the chromatic polynomial is {\em unimodal}: it goes up, hits a peak, and then goes down.  
Read's conjecture was proved by June Huh in a 2012 paper \cite{Huh_JAMS} making heavy use of methods from algebraic geometry.
Huh's result was subsequently refined and generalized by Huh and Katz \cite{HK}, again using substantial doses of algebraic geometry.
Both papers in fact establish log-concavity of the coefficients, which is stronger than unimodality.  

The breakthroughs of Huh and Huh--Katz left open the more general Rota--Welsh conjecture where graphs are generalized to (not necessarily representable) matroids and the chromatic polynomial of a graph is replaced by the characteristic polynomial of a matroid.   The Huh and Huh--Katz techniques are not applicable in this level of generality, since there is no underlying algebraic geometry to which to relate the problem.  But in 2015 Adiprasito, Huh, and Katz \cite{AHK} announced a proof of the Rota--Welsh conjecture based on a novel approach motivated by but not making use of any results from algebraic geometry. The authors first prove that the Rota--Welsh conjecture would follow from combinatorial analogues of the Hard Lefschetz Theorem and Hodge-Riemann relations in algebraic geometry.  
They then implement an elaborate inductive procedure to prove the combinatorial Hard Lefschetz Theorem and Hodge-Riemann relations using purely combinatorial arguments.

We will survey these developments.
\end{abstract}

\maketitle


\section{Unimodality and Log-Concavity}

A sequence $a_0,\ldots,a_d$ of real numbers is called {\em unimodal} if there is an index $i$ 
such that
\[
a_0 \leq \cdots \leq a_{i-1} \leq a_i \geq a_{i+1} \geq \cdots \geq a_d.
\]

There are numerous naturally-occurring unimodal sequences in algebra, combinatorics, and geometry.  For example:

\begin{example} (Binomial coefficients) \label{ex:binomial}
The sequence of binomial coefficients $\binom{n}{k}$ for $n$ fixed and $k=0,\ldots,n$ (the $n^{\rm th}$ row of Pascal's triangle) is unimodal.
\end{example}

The sequence $\binom{n}{0},\binom{n}{1},\ldots,\binom{n}{n}$ has a property which is in fact stronger than unimodality: it is {\em log-concave}, meaning
that $a_i^2 \geq a_{i-1} a_{i+1}$ for all $i$.   Indeed,
\[
\frac{\binom{n}{k}^2}{\binom{n}{k-1} \binom{n}{k+1}} = \frac{(k+1)(n-k+1)}{k(n-k)} > 1.
\]

It is a simple exercise to prove that a log-concave sequence of positive numbers is unimodal.

Some less trivial, but still classical and elementary, examples of log-concave (and hence unimodal) sequences are the Stirling numbers of the first and second kind.  

\begin{example} (Stirling numbers) \label{ex:stirling}
The {\em Stirling numbers of the first kind}, denoted $s(n,k)$, are the coefficients which appear when one writes falling factorials $x\downarrow_n =  x(x-1)\cdots(x-n+1)$ as polynomials in $x$:
\[
x\downarrow_n = \sum_{k=0}^n s(n,k) x^k.
\]
This sequence of integers alternates in sign.  The {\em signless Stirling numbers of the first kind} $s^+(n,k) = |s(n,k)| = (-1)^{n-k} s(n,k)$ enumerate the number of permutations of $n$ elements having exactly $k$ disjoint cycles.

The {\em Stirling numbers of the second kind}, denoted $S(n,k)$, invert the Stirling numbers of the first kind in the sense that
\[
\sum_{k=0}^n S(n,k) (x)_k = x^n.
\]
Their combinatorial interpretation is that $S(n,k)$ counts the number of ways to partition an $n$ element set into $k$ non-empty subsets.

For fixed $n$ (with $k$ varying from $0$ to $n$), both $s^+(n,k)$ and $S(n,k)$ are log-concave and hence unimodal.
\end{example}

Another example, proved much more recently through a decidedly less elementary proof, concerns the sequence of coefficients of the chromatic polynomial of a graph.  This example will be the main focus of our paper.

\begin{example}(Coefficients of the chromatic polynomial) \label{ex:chromatic}
Let $G$ be a finite graph\footnote{We allow loops and parallel edges.}.  In 1912, George Birkhoff defined $\chi_G(t)$ to be the number of proper colorings of $G$ using $t$ colors (i.e., the number of functions $f : V(G) \to \{ 1,\ldots,t \}$ such that $f(v) \neq f(w)$ whenever
$v$ and $w$ are adjacent in $G$), and proved that $\chi_G(t)$ is a {\em polynomial} in $t$, called the {\bf chromatic polynomial} of $G$.  

For example,\footnote{Looking at the analogy between the formulas $\chi_T$ and $\chi_{K_n}$, and between $s(n,k)$ and $S(n,k)$, it may be reasonable to think of $(-1)^{n-k} \binom{n}{k}$ as a ``binomial coefficient of the first kind'' and of the usual binomial coefficients as being of the ``second kind''.  This fits in neatly with the ``inversion formulas'' $(1+x)^n = \sum \binom{n}{k} x^k$ and $x^n = \sum (-1)^{n-k} \binom{n}{k} (1+x)^k$.} if $G=T$ is a tree on $n$ vertices then the chromatic polynomial of $G$ is
\[
\chi_T(t) = t(t-1)^{n-1} = \sum_{k=1}^n (-1)^{n-k} \binom{n-1}{k-1} t^k.
\]

If $G=K_n$ is the complete graph on $n$ vertices, then 
\[
\chi_{K_n}(t) = t(t-1)\cdots (t-n+1) = \sum_{k=1}^n s(n,k) t^k.
\]

And if $G$ is the Petersen graph, depicted in Figure~\ref{fig:Petersen}, then 
\[
\chi_G(t) = t^{10}-15 t^9+105 t^8-455 t^7+1353 t^6-2861 t^5+4275 t^4-4305 t^3+2606 t^2-704 t.
\]

\begin{figure}
\centering
\includegraphics[scale=0.4]{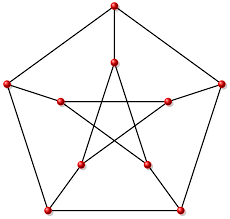}
\caption{The Petersen graph}
\label{fig:Petersen}
\end{figure}

Ronald Reed conjectured in 1968 that for any graph $G$ the (absolute values of the) coefficients of $\chi_G(t)$ form a unimodal sequence, and a few years letter Stuart Hoggar conjectured that the coefficients in fact form a log-concave sequence\footnote{Log-concavity implies unimodality for the coefficients of $\chi_G(t)$ by the theorem of Rota mentioned at the end of \S\ref{sec:MobiusFunction}.}.
Both conjectures were proved only relatively recently by June Huh \cite{Huh_JAMS}.
\end{example}

Another interesting and relevant example concerns linearly independent sets of vectors:

\begin{example} \label{ex:welsh}
Let $k$ be a field, let $V$ be a vector space over $k$, and let $A$ be a finite subset of $V$.
Dominic Welsh conjectured that $f_i(A)$ is a log-concave sequence, where $f_i(A)$ is the number of linearly independent subsets of $A$ of size $i$.
For example, if $k = {\mathbb F}_2$ is the field of $2$ elements, $V = {\mathbb F}_2^3$, and $A = V \backslash \{ 0 \}$, then
\[
f_0(A)=1, f_1(A)=7, f_2(A)=21, f_3(A)=28.
\]
This conjecture 
is a consequence of the recent work of Huh--Katz \cite{HK} (cf.~\cite{Lenz}).
\end{example}

Finally, we mention an example of an apparently much different nature coming from algebraic geometry:

\begin{example} (Hard Lefschetz Theorem) \label{ex:hardlefschetz}
Let $X$ be an irreducible smooth projective algebraic variety of dimension $n$ over the field $\CC$ of complex numbers, and let $\beta_i = {\rm dim \;} H^i(X,\CC)$ be the {\em $i^{\rm th}$ Betti number} of $X$. (Here $H^*(X,\CC)$ denotes the singular cohomology groups of $X$.)  Then the two sequences $\beta_0,\beta_2,\ldots,\beta_{2n}$ and $\beta_1,\beta_3,\ldots,\beta_{2n-1}$ are symmetric and unimodal.
Moreover, this remains true if we replace the hypothesis that $X$ is smooth by the weaker hypothesis that $X$ has only finite quotient singularities, meaning that $X$ looks locally (in the analytic topology) like the quotient of $\CC^n$ by a finite group of linear transformations.

The symmetry of the $\beta_i$'s is a classical result in topology known as {\em Poincar{\'e} duality}.
 And one has the following important strengthening (given symmetry) of unimodality: there is an element $\omega \in H^2(X,\CC)$ such that for $0 \leq i \leq n$, multiplication by $\omega^{n-i}$ defines an isomorphism from $H^i(X,\CC)$ to $H^{2n-i}(X,\CC)$.
This result is called the {\em Hard Lefschetz Theorem}.  In the smooth case, it is due to Hodge; for varieties with finite quotient singularities, it is due to Saito and uses the theory of perverse sheaves.

For varieties $X$ with arbitrary singularities, the Hard Lefschetz Theorem still holds if one replaces singular cohomology by the {\em intersection cohomology} of Goresky and MacPherson (cf. \cite{CDM_Bulletin}).
\end{example}

Surprisingly, all five of the above examples are in fact related.  We have already seen that Example~\ref{ex:binomial}, as well as Example~\ref{ex:stirling} in the case of Stirling numbers of the first kind, are special cases of Example~\ref{ex:chromatic}.
We will see in the next section that Examples~\ref{ex:chromatic} and \ref{ex:welsh} both follow from a more general result concerning {\bf matroids}.  And the proof of this theorem about matroids will involve, as one of its key ingredients, a combinatorial analogue of the Hard Lefschetz Theorem (as well as the Hodge-Riemann relations, about which we will say more later).

\section{Matroids}

Our primary references for this section are \cite{Oxley} and \cite{Welsh}.

\subsection{Independence axioms}
Matroids were introduced by Hassler Whitney as a combinatorial abstraction of the notion of linear independence of vectors.
There are many different (``cryptomorphic'') ways to present the axioms for matroids, all of which turn out to be non-obviously equivalent to one another. 
For example, instead of using linear independence one can also define matroids by abstracting the notion of {\em span}.
We will give a brief utilitarian introduction to matroids, starting with the independence axioms.

\begin{definition}
\label{def:independence}
(Independence Axioms) A {\bf matroid} $M$ is a finite set $E$ together with a collection ${\mathcal I}$ of subsets of $E$, called the {\em independent sets} of the matroid, such that:
\begin{itemize}
\item[(I1)] The empty set is independent.
\item[(I2)] Every subset of an independent set is independent.
\item[(I3)] If $I,J$ are independent sets with $|I| < |J|$, then there exists $y \in J \backslash I$ such that $I \cup \{ y \}$ is independent.
\end{itemize}
\end{definition}

\subsection{Examples}

\begin{example} \label{ex:linear} (Linear matroids) 
Let $V$ be a vector space over a field $k$, and let $E$ be a finite subset of $V$.  Define ${\mathcal I}$ to be the collection of linearly independent subsets of $E$.
Then ${\mathcal I}$ satisfies (I1)-(I3) and therefore defines a matroid.  
Slightly more generally (because we allow repetitions), if $E = \{ 1,\ldots, m \}$ and A is an $n \times m$ matrix with entries in $k$,
a subset of $E$ is called independent iff the corresponding columns of $A$ are linearly independent over $k$.  
We denote this matroid by $M_k(A)$. 
A matroid of the form $M_k(A)$ for some $A$ is called {\bf representable} over $k$.

By a recent theorem of Peter Nelson \cite{Nelson}, asymptotically 100\% of all matroids are not representable over any field.
\end{example}

\begin{example} (Graphic matroids)
Let $G$ be a finite graph, let $E$ be the set of edges of $G$, and let ${\mathcal I}$ be the collection of all subsets of $E$ which do not contain a cycle. Then ${\mathcal I}$ satisfies (I1)-(I3) and hence defines a matroid $M(G)$.
The matroid $M(G)$ is {\bf regular}, meaning that it is representable over every field $k$.  By a theorem of Whitney, if $G$ is 3-connected (meaning that $G$ remains connected after removing any two vertices) then $M(G)$ determines the isomorphism class of $G$.
\end{example}

\begin{example}  \label{ex:uniform} (Uniform matroids)
Let $E = \{ 1,\ldots, m \}$ and let $r$ be a positive integer.  The {\bf uniform matroid} $U_{r,m}$ is the matroid on $E$ whose independent sets are the subsets of $E$ of cardinality at most $r$.  For each $r,m$ there exists $N=N(r,m)$ such that $U_{r,m}$ is representable over every field having at least $N(r,m)$ elements.
\end{example}

\begin{figure}
\centering
\includegraphics[scale=0.3]{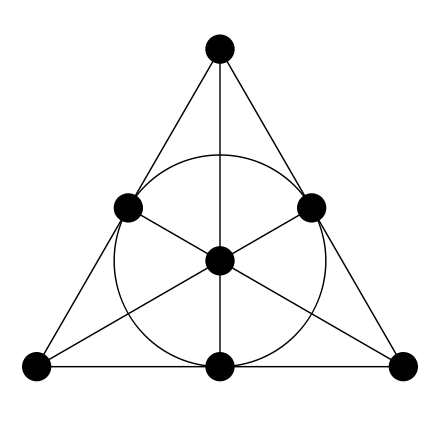}
\caption{The Fano matroid}
\label{fig:Fano}
\end{figure}

\begin{example} \label{ex:Fano} (Fano matroid) 
Let $E={\mathbb P}^2({\mathbf F}_2)$ be the projective plane over the 2-element field; the seven elements of $E$ can be identified with the dots in
Figure~\ref{fig:Fano}.

Define ${\mathcal I}$ to be the collection of subsets of $E$ of size at most 3 which are not one of the 7 lines in ${\mathbb P}^2({\mathbf F}_2)$ (depicted as six straight lines and a circle in Figure~\ref{fig:Fano}). 
Then ${\mathcal I}$ satisfies (I1)-(I3) and determines a matroid called the {\em Fano matroid}.  This matroid is representable over ${\mathbf F}_2$ but not over any field of characteristic different from 2. In particular, the Fano matroid is not graphic.
\end{example}

\begin{figure}
\centering
\includegraphics[scale=0.2]{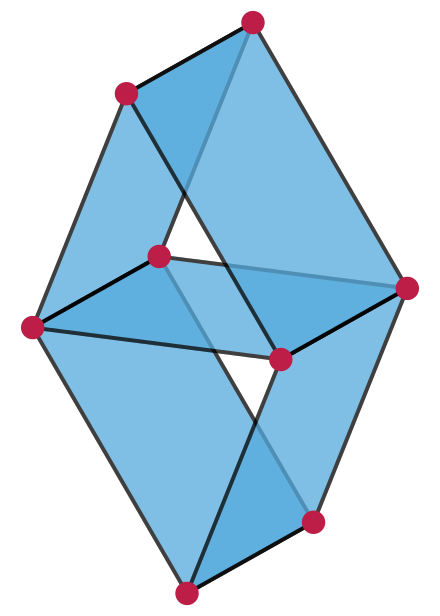}
\caption{The Vamos matroid}
\label{fig:Vamos}
\end{figure}

\begin{example} (Vamos matroid) 
Let $E$ be the 8 vertices of the cuboid shown in Figure~\ref{fig:Vamos}.
Define ${\mathcal I}$ to be the collection of subsets of $E$ of size at most 4 which are not one of the five square faces in the picture.
Then ${\mathcal I}$ satisfies (I1)-(I3) and determines a matroid called the {\em Vamos matroid} which is not representable over any field.
\end{example}

\subsection{Circuits, bases, and rank functions}

A subset of $E$ which is not independent is called dependent.  A minimal dependent set is called a {\bf circuit}, and a maximal independent set is called a {\bf basis}.
As in linear algebra, all bases of $M$ have the same cardinality; this number is called the {\bf rank} of the matroid $M$, and is denoted $r(M)$.
More generally, if $A$ is a subset of $E$, we define the {\bf rank} of $A$, denoted $r_M(A)$ or just $r(A)$,
to be the maximal size of an independent subset of $A$.

One can give cryptomorphic axiomatizations of matroids in terms of circuits, bases, and rank functions.
For the sake of brevity we refer the interested reader to \cite{Oxley}.

\subsection{Duality}
If $M = (E,{\mathcal I})$ is a matroid, let ${\mathcal I}^*$ be the collection of subsets $A\subseteq E$ such that $E \backslash A$ contains a basis $B$ for $M$. It turns out that ${\mathcal I}^*$ satisfies axioms (I1),(I2), and (I3) and thus 
$M^* = (E,{\mathcal I}^*)$ is a matroid, called the {\bf dual matroid} of $M$.


If $M=M(G)$ is the matroid associated to a planar graph $G$, then $M^*$ is the matroid associated to the {\em planar dual} of $G$. A theorem of Whitney asserts, conversely, that if $G$ is a connected graph for which the dual matroid $M(G)^*$ is graphic, then $G$ is planar.

\subsection{Deletion and Contraction}

Given a matroid $M$ on $E$ and $e \in E$, we write $M \backslash e$ for the matroid on $E \backslash \{ e \}$ 
whose independent sets are the independent sets of $M$ not containing $e$.

We write $M / e$ for the matroid on $E \backslash \{ e \}$ 
such that $I$ is independent in $M / e$ if and only if $I = J \backslash \{ e \}$ with $J$ independent in $M$ and $e \in J$.

We call these operations on matroids {\bf deletion} and {\bf contraction}, respectively.
Deletion and contraction are dual operations, in the sense that $(M \backslash e)^* = M^* / e$ and $(M / e)^* = M^* \backslash e$.

If $M$ is a graphic matroid, deletion and contraction correspond to the usual notions in graph theory.

\subsection{Spans}
We defined matroids in terms of independent sets, which abstract the notion of linear independence.
We now focus on a different way to define / characterize matroids in terms of {\em closure operators}, which abstract the notion of {\em span} in linear algebra.

Let $2^E$ denote the power set of $E$.

\begin{definition}
\label{def:closure}
(Span Axioms) A {\bf matroid} $M$ is a finite set $E$ together with a function ${\rm cl}: 2^E \to 2^E$ such that for all $X,Y \subseteq E$ and $x,y \in E$:
\begin{itemize}
\item[(S1)] $X \subseteq {\rm cl}(X)$.
\item[(S2)] If $Y \subseteq X$ then ${\rm cl}(Y) \subseteq {\rm cl}(X)$.
\item[(S3)] ${\rm cl}({\rm cl}(X)) = {\rm cl}(X)$.
\item[(S4)] If $y \in {\rm cl}(X \cup \{ x \})$ but $y \not\in {\rm cl}(X)$, then $x \in {\rm cl}(X \cup \{ y \})$.
\end{itemize}
\end{definition}

For example, if $M$ is a linear matroid as in Example~\ref{ex:linear} then ${\rm cl}(X)$ is just the span of $X$ in $V$.

The {\em exchange axiom} (S4) captures our intuition of a ``geometry'' as a collection of incidence relations
\[
\{ {\rm point} \} \subset \{ {\rm line} \} \subset \{ {\rm plane} \} \subset \cdots
\]

For example, if $L$ is a line in an $r$-dimensional projective space $\PP^{r}_k$ over a field $k$ and $p,q \in \PP^{r}_k \backslash L$,
then $q$ lies in the span of $L \cup \{ p \} \iff$ $p$ lies in the span of $L \cup \{ q \}$ $\iff$ $p,q,L$ are {\em coplanar}.

The relation between Definitions~\ref{def:independence} and \ref{def:closure} is simple to describe: given a matroid in the sense of Definition~\ref{def:independence}, we define ${\rm cl}(X)$ to be the set of all $x \in E$ such that $r(X \cup \{ x \}) = r(X)$. 
Conversely, given a matroid in the sense of Definition~\ref{def:closure}, we define a subset $I$ of $E$ to be independent if and only if $x \in I$ implies $x \not\in {\rm cl}(I \backslash \{ x \})$.

A subset $X$ of $E$ is said to {\bf span} $M$ if ${\rm cl}(X) = E$.  As in the familiar case of linear algebra, one can show in general that $X$ is a basis (i.e., a maximal independent set) if and only if $X$ is independent and spans $E$.

\subsection{Flats}

A subset $X$ of $E$ is called a {\bf flat} (or a {\em closed subset}) if $X = {\rm cl}(X)$.  

\begin{example} \label{ex:linearflats} (Linear matroids) 
Let $V$ be a vector space and let $E$ be a finite subset of $V$.  
A subset $F$ of $E$ is a flat of the corresponding linear matroid if and only if there is no vector in $E \backslash F$ contained in the linear span of $F$.

Alternatively, let $M = M_k(A)$ be represented by an $r \times m$ matrix $A$ of rank $r$ with entries in $k$, and let $V \subseteq k^m$ be the row space of $A$.  
Let $E = \{ 1,\ldots, m \}$, and for $I \subseteq E$ let $L_I$ be the ``coordinate flat''
\[
L_I = \{ x=(x_1,\ldots,x_m) \in k^m \; : \; x_i = 0 \; {\rm for} \; i \in I \}.
\]
Then for $I \subseteq E$ we have $r_M(I) = {\rm dim}(V) - {\rm dim}(V \cap L_I)$, and $I$ is a flat of $M$ if and only if $V \cap L_J \subsetneq V \cap L_I$ for all $J \supsetneq I$.
In particular, $V \cap L_I = V \cap L_F$, where $F$ is the smallest flat of $M$ containing $I$.
\end{example}

\begin{example} (Graphic matroids)
Let $G$ be a connected finite graph, and let $M(G)$ be the associated matroid.
Then a subset $F$ of $E$ is a flat of $M(G)$ if and only if there is no edge in $E \backslash F$ whose endpoints are connected by a path in $F$.
\end{example}

\begin{example} (Fano matroid) 
In the Fano matroid, the flats are $\emptyset,E$, and each of the 7 points and 7 lines in Figure~\ref{fig:Fano}.
\end{example}

Every maximal chain of flats of a matroid $M$ has the same length, which coincides with the rank of $M$.

One can give a cryptomorphic axiomatization of matroids in terms of flats.
To state it, we say that a flat $F'$ {\bf covers} a flat $F$ if $F \subsetneq F'$ and there are no intermediate flats between $F$ and $F'$.

\begin{definition}
\label{def:flats}
(Flat Axioms) A {\bf matroid} $M$ is a finite set $E$ together with a collection of subsets of $E$, called {\bf flats}, such that:\footnote{For the geometric intuition behind axiom (F3), note that given a line in $\RR^3$, the planes which contain this line (minus the line itself) partition the remainder of $\RR^3$.}
\begin{itemize}
\item[(F1)] $E$ is a flat.
\item[(F2)] The intersection of two flats is a flat.
\item[(F3)] If $F$ is a flat and $\{ F_1,F_2,\ldots,F_k \}$ is the set of flats that cover $F$, then $\{ F_1 \backslash F, F_2 \backslash F, \ldots, F_k \backslash F \}$ partitions $E \backslash F$.
\end{itemize}
\end{definition}

We have already seen how to define flats in terms of a closure operator.  To go the other way, one defines the closure of a set $X$ to be the intersection of all flats containing $X$.

\subsection{Simple matroids}

A matroid $M$ is called {\bf simple} if every dependent set has size at least $3$.  Equivalently, a matroid is simple if and only if it has no:

\begin{itemize}
\item {\bf Loops} (elements $e \in {\rm cl}(\emptyset)$); or
\item {\bf Parallel elements}\footnote{The terms ``loop'' and ``parallel element'' come from graph theory.} (elements $e,e'$ with $e' \in {\rm cl}(e)$).
\end{itemize}

Every matroid $M$ has a canonical {\bf simplification} $\hat{M}$ obtained by removing all loops and identifying parallel elements (with the obvious resulting notions of independence, closure, etc.).
A simple matroid is also called a {\bf combinatorial geometry}.

For future reference, we define a {\bf coloop} of a matroid $M$ to be a loop of $M^*$.  Equivalently, a loop is an element $e \in E$ which does not belong to any basis of $M$, and a coloop is an element $e \in E$ which belongs to every basis of $M$.

\subsection{The Bergman fan of a matroid} \label{sec:BergmanFan}

Let $E=\{ 0,1,\ldots,n \}$ and let $M$ be a matroid on $E$.
The Bergman fan of $M$ is a certain collection of cones in the $n$-dimensional Euclidean space $N_{\RR} = \RR^E / \RR (1,1,\ldots,1)$ which 
carries the same combinatorial information as $M$.
Bergman fans show up naturally in the context of tropical geometry, where they are also known (in the ``trivially-valued case'') as {\em tropical linear spaces}.

For $S \subseteq E$, let $e_S = \sum_{i \in S} e_i \in N_\RR$, where $e_i$ is the basis vector of $\RR^E$ corresponding to $i$.
Note that $e_E = 0$ by the definition of $N_\RR$.
Let $F_\bullet = \{ F_1 \subsetneq F_2 \subsetneq \cdots \subsetneq F_k \}$ be a $k$-step flag of non-empty proper flats of $M$.
We define the corresponding cone $\sigma_{F_\bullet} \subseteq N_\RR$ to be the nonnegative span of the $e_{F_i}$ for $i=1,\ldots,k$.

\begin{definition}
The {\bf Bergman fan} $\Sigma_M$ of $M$ is the collection of cones $\sigma_{F_\bullet}$ as $F_\bullet$ ranges over all flags of non-empty proper flats of $M$.
\end{definition}

\begin{example} \label{ex:tropicalline}
The Bergman fan of the uniform matroid $U_{2,3}$ has a zero-dimensional cone given by the origin in $\RR^2$ and three $1$-dimensional cones given by rays from the origin in the directions of $\bar{e}_1 = (1,0), \bar{e}_2 = (0,1)$, and $\bar{e}_3 = -(\bar{e}_1 + \bar{e}_2) = (-1,-1)$.  This is the well-known ``tropical line'' in $\RR^2$ with vertex at the origin, see Figure~\ref{fig:tropicalline}.
\begin{figure}
\centering
\includegraphics[scale=0.4]{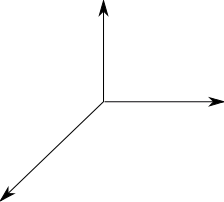}
\caption{The Bergman fan of $U_{2,3}$}
\label{fig:tropicalline}
\end{figure}
\end{example}

\begin{example} \label{ex:permutohedron}
Let $U = U_{n+1,n+1}$ be the rank $n+1$ uniform matroid on $E=\{ 0,1,\ldots,n \}$.
Every subset of $E$ is a flat, so the top-dimensional cones of $\Sigma_U$ are the nonnegative spans of 
\[
\{ e_{i_0}, e_{i_0} + e_{i_1}, \ldots, e_{i_0} + e_{i_1} + \cdots + e_{i_{n-1}} \}
\]
for every permutation $i_0,\ldots,i_n$ of $0,\ldots,n$.
The fan $\Sigma_U$ is the normal fan to the {\bf permutohedron} $P_n$, which by definition is the convex hull of $(i_0,\ldots,i_n)$ over all permuations $i_0,\ldots,i_n$ of $0,\ldots,n$, viewed as a polytope in the dual vector space to $N_\RR$, see Figure~\ref{fig:permutohedron}.)
\begin{figure}
\centering
\includegraphics[scale=0.3]{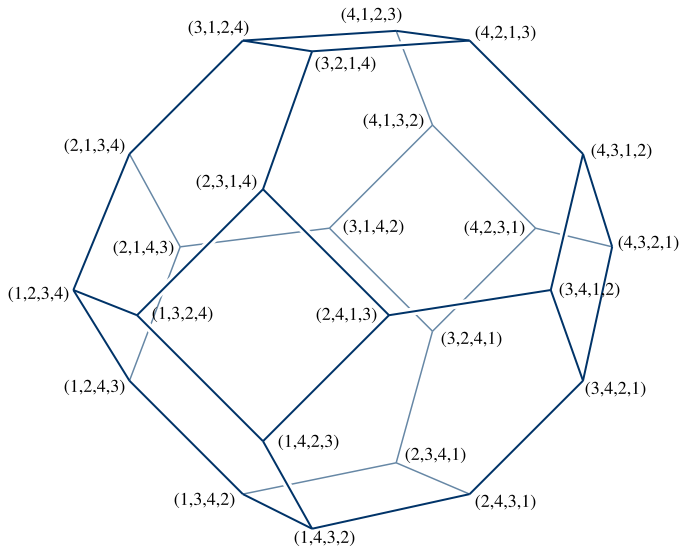}
\caption{The 3-dimensional permutohedron.  Note that here we take $E=\{ 1,2,3,4 \}$, instead of $\{ 0,1,2,3 \}$ as in Example~\ref{ex:permutohedron}, but this is immaterial since we work in $\RR^E / \RR (1,1,\ldots,1)$.}
\label{fig:permutohedron}
\end{figure}
The fan $\Sigma_U$ plays a central and recurring role in \cite{AHK}.
For further information on some of the remarkable combinatorial and geometric properties of $P_n$, see \cite{ZieglerPolytopes}.
\end{example}

The isomorphism class of matroid $M$ determines and is determined by its Bergman fan.\footnote{One can in fact give a cryptomorphic characterization of matroids via their Bergman fans, using the flat axioms (F1)-(F3): a rational polyhedral fan $\Sigma$ in $N_\RR$ is the Bergman fan of a matroid on $E$ if and only if it is {\em balanced} and has {\em degree one} as a tropical cycle \cite{Fink}.}

\section{Geometric lattices and the characteristic polynomial}

In this section we define the characteristic polynomial of a matroid $M$ in terms of the lattice of flats of $M$.
Our primary references are \cite{Welsh} and \cite[Chapters 7-8]{WhiteCG}.

\subsection{Geometric lattices}

The set ${\mathcal L}(M)$ of flats of a matroid $M$ together with the inclusion relation forms a {\bf lattice}, i.e., a partially ordered set in which every two elements $x,y$ have both a {\bf meet} (greatest lower bound) $x \wedge y$ and a {\bf join} 
(least upper bound) $x \vee y$.  Indeed, if $X$ and $Y$ are flats then we can define $X \wedge Y$ as the intersection of $X$ and $Y$ and $X \vee Y$ as the closure of the union of $X$ and $Y$.

\begin{example}
Flats of the uniform matroid $U_{n+1,n+1}$ can be identified with subsets of $\{ 0,1,\ldots,n \}$, and with this identification the lattice of flats of $U_{n+1,n+1}$ is the {\bf Boolean lattice} $B_{n+1}$ consisting of subsets of $\{ 0,1,\ldots,n \}$ partially ordered by inclusion.
\end{example}

\begin{example}
Flats of the complete graph $K_n$ can be identified with partitions of $\{ 1,\ldots, n \}$, and with this identification the lattice of flats of the graphic matroid $M(K_n)$ is isomorphic to the {\bf partition lattice} $\Pi_n$ consisting of partitions of  $\{ 1,\ldots, n \}$ partially ordered by refinement.
\end{example}

If $L$ is a lattice and $x,y \in L$, we say that $y$ {\bf covers} $x$ if
$x < y$ and whenever $x \leq z \leq y$ we have either $z = x$ or $z=y$.  A finite lattice has a minimal element $0_L$ and a maximal element $1_L$.  An {\bf atom} is an element which covers $0_L$.

The lattice of flats $L = {\mathcal L}(M)$ has the following properties:

\begin{itemize}
\item[(L1)] $L$ is {\bf semimodular}, i.e., if $x,y \in L$ both cover $x \wedge y$ then $x \vee y$ covers both $x$ and $y$.
\item[(L2)] $L$ is {\bf atomic}, i.e., every $x \in L$ is a join of atoms.
\end{itemize}

A lattice satisfying (L1) and (L2) is called a {\bf geometric lattice}.  By a theorem of Garrett Birkhoff (the son of George), every geometric lattice is of the form ${\mathcal L}(M)$ for some matroid $M$.  However, the matroid $M$ is not unique, because if $\hat{M}$ is the simplification of $M$ then ${\mathcal L}(M) = {\mathcal L}(\hat{M})$.  Birkhoff proves that this is in fact the only ambiguity, i.e., the map $M \mapsto {\mathcal L}(M)$ gives a bijection between isomorphism classes of simple matroids and isomorphism classes of geometric lattices.
Thus, at least up to simplification, (L1) and (L2) give another cryptomorphic characterization of matroids.

If $F$ is a flat of a matroid $M$, the maximal length $\ell$ of a chain $F_0 \subset F_1 \subset \cdots \subset F_\ell = F$ of flats coincides with the rank $r_M(F)$ of $F$.
This allows us to define the rank function on $M$, restricted to the set of flats, purely in terms of the lattice ${\mathcal L}(M)$.
We write $r_L$ for the corresponding function on an arbitrary geometric lattice $L$.

\subsection{The M{\"o}bius function of a poset} \label{sec:MobiusFunction}

There is a far-reaching combinatorial result known as
the M{\"o}bius Inversion Formula which holds in an arbitrary finite poset $P$.
It simultaneously generalizes, among other things, the Inclusion-Exclusion Principle, the usual number-theoretic M{\"o}bius Inversion Formula, and the Fundamental Theorem of Difference Calculus.

There is a unique function $\mu_P : P \times P \to \ZZ$, called the {\bf M{\"o}bius function} of $P$, satisfying $\mu_P(x,x)=1$, $\mu_P(x,y)=0$ if $x \not\leq y$, and 
\[
\sum_{x \leq z \leq y} \mu_P(x,z) = 0.
\]
if $x < y$.  
Note that $\mu_P(x,y)=-1$ if $y$ covers $x$.

The {\bf M{\"o}bius Inversion Formula} states that if $f$ is a function from a finite poset $P$ to an abelian group $H$, and if we define $g(y) = \sum_{x \leq y} f(x)$ for all $y \in P$, then
\[
f(y) = \sum_{x \leq y} \mu_P(x,y) g(x).
\]

If $P=L$ is a finite {\em lattice}, the M{\"o}bius function satisfies {\bf Weisner's theorem}, which gives a ``shortcut'' for the recurrence defining $\mu$: if $0_L \neq x \in L$ then 
\[
\sum_{y \in L \; : \; x \vee y = 1_L } \mu_L(0_L,y) = 0.
\]

If $L$ is moreover a {\em geometric lattice}, it is a theorem of Rota that the M{\"o}bius function of $L$ is non-zero and alternates in sign.  More precisely, if $x \leq y$ in $L$ then
\[
(-1)^{r_L(y) - r_L(x)} \mu_L(x,y) > 0.
\]

\subsection{The characteristic polynomial} \label{sec:charpoly}

The chromatic polynomial of a graph $G$ satisfies the {\em deletion-contraction relation}:
\[
\chi_G(t) = \chi_{G \backslash e}(t) - \chi_{G / e}(t).
\]

Indeed, the equivalent formula $\chi_{G \backslash e}(t) =  \chi_G(t) + \chi_{G / e}(t)$ just says that 
the proper colorings of $G \backslash e$ can be partitioned into those where the endpoints of $e$ are colored differently (giving a proper coloring of $G$) or the same (giving a proper coloring of $G/e$).


This formula is not only useful for calculating $\chi_G(t)$, it is also the simplest way to prove that $\chi_G(t)$ is a polynomial in $t$
(by induction on the number of edges).
In addition, this formula for $\chi_G(t)$ suggests an extension to arbitrary matroids.
This can be made to work, but it is not obvious that this recursive procedure is always well-defined.  
So it is more convenient to proceed as follows. 

First, note that the chromatic polynomial of a graph $G$ is identically zero by definition if $G$ has a loop edge.  So we will define $\chi_M(t)=0$ for any
matroid with a loop.  We may thus concentrate on loopless matroids.  Note that a matroid $M$ is loopless if and only if $\emptyset$ is a flat of $M$.

\begin{definition} \label{def:charpoly}
Let $M$ be a loopless matroid with lattice of flats $L$.  The {\bf characteristic polynomial} of $M$ is
\begin{equation}  \label{eq:chi1}
\chi_M(t) = \sum_{F \in L} \mu_L(\emptyset,F) t^{r(M) - r(F)}.
\end{equation}
\end{definition}

In particular, if $M$ is loopless then $\chi_M(t) = \chi_{\hat{M}}(t)$, where $\hat{M}$ denotes the simplification of $M$.  

The motivation behind (\ref{eq:chi1}) may be unclear to the reader at this point.  
In the representable case, at least, there is a ``motivic'' interpretation of (\ref{eq:chi1}) which some will find illuminating; see \S\ref{sec:motivic}.

There is also a (simpler-looking but sometimes not as useful) expression for $\chi_M(t)$ in terms of a sum over {\em all} subsets of $E$, not just flats.

\begin{proposition} \label{prop:simplechi}
If $M$ is any matroid,
\begin{equation}  \label{eq:chi2}
\chi_M(t) = \sum_{A \subseteq E} (-1)^{|A|} t^{r(M) - r(A)}.
\end{equation}
\end{proposition}

If $M_1,M_2$ are matroids on $E_1$ and $E_2$, respectively, and $E_1 \cap E_2 = \emptyset$, we define the {\bf direct sum} $M_1 \oplus M_2$ to be the matroid on $E_1 \cup E_2$ whose flats are all sets of the form $F_1 \cup F_2$ where $F_i$ is a flat of $M_i$ for $i=1,2$.
The following result gives an important characterization of the characteristic polynomial.

\begin{theorem} \label{thm:charpolychar}
Let $M$ be a matroid.
\begin{itemize}
\item[$(\chi 1)$] If $e$ is neither a loop nor a coloop of $M$, then $\chi_M(t) =  \chi_{M \backslash e}(t) - \chi_{M / e}(t)$.
\item[$(\chi 2)$] If $M=M_1 \oplus M_2$ then $\chi_M(t) = \chi_{M_1}(t) \chi_{M_2}(t)$.
\item[$(\chi 3)$] If $M$ contains a loop then $\chi_M(t)=0$, and if $M$ consists of a single coloop then $\chi_M(t) = t-1$.
\end{itemize}
Furthermore, the characteristic polynomial is the {\em unique} function from matroids to integer polynomials satisfying $(\chi 1)$-$(\chi 3)$.
\end{theorem}

In particular, it follows from Theorem~\ref{thm:charpolychar} that if $G$ is a graph then the chromatic polynomial $\chi_G(t)$ of $G$ satisfies $\chi_G(t) = t^{c(G)} \chi_{M(G)}(t)$, where $c(G)$ is the number of connected components of $G$.
(The extra factor of $t$ when $G$ is connected comes from the fact that the graph with two vertices and one edge has chromatic polynomial $t(t-1)$, whereas the corresponding matroid, which consists of a single coloop, has characteristic polynomial $t-1$.
Note that since no graph can be $0$-colored, $\chi_G(0)=0$ for every graph $G$ and hence the chromatic polynomial is always divisible by $t$.)

The characteristic polynomial of $M$ is monic of degree $r = r(M)$, so we can write 
\[
\chi_M(t) = w_0(M)t^{r} + w_1(M) t^{r - 1} + \cdots + w_{r}(M)
\]
with $w_0(M) = 1$ and $w_k(M) \in \ZZ$.
By Rota's theorem, the coefficients of $\chi_M(t)$ alternate in sign, i.e.,
\[
w_k(M)^+ := |w_k(M)| = (-1)^k w_k(M).
\]

The numbers $w_k(M)$ (resp. $w_k^+(M)$) are called the {\bf Whitney numbers of the first kind} (resp. {\bf unsigned Whitney numbers of the first kind}) for $M$.
The recent work of Adiprasito--Huh--Katz \cite{AHK} establishes:

\begin{theorem} \label{thm:logconcavefirstkind}
For any matroid $M$, the unsigned Whitney numbers of the first kind $w_k^+(M)$ form a log-concave sequence.
\end{theorem}

Note that it is enough to prove the theorem for simple matroids, i.e., combinatorial geometries, since the characteristic polynomial of a loopless matroid equals that of its simplification.

Actually, Adiprasito, Huh, and Katz study the so-called {\bf reduced characteristic polynomial} of $M$.
If $|E| \geq 1$ then $\chi_M(1)=0$ (e.g., if $G$ is a graph with at least one edge then $G$ has no proper one-coloring!).
Thus we may write $\chi_M(t) = (t-1) \bar{\chi}_M(t)$ with $\bar{\chi}_M(t) \in \ZZ[t]$.
The reduced characteristic polynomial $\bar{\chi}_M(t)$ is the ``projective analogue'' of $\chi_M(t)$
(cf.~\S\ref{sec:motivic} below).
It is an elementary fact that log-concavity of the (absolute values of the) coefficients of $\bar{\chi}_M(t)$ implies log-concavity for $\chi_M(t)$.
So in order to prove Theorem~\ref{thm:logconcavefirstkind} one can replace the $w_k^+(M)$ by their projective analogues $m_k(M)$.

\subsection{Tutte-Grothendieck invariants}
Our primary reference for this section and \S\ref{sec:motivic} is \cite{KatzMatroidAG}.

One can generalize the characteristic polynomial of a matroid by relaxing the condition that it vanishes on matroids containing loops.

The {\bf Tutte polynomial} of a matroid $M$ on $E$ is the two-variable polynomial
\[
T_M(x,y) = \sum_{A \subseteq E} (x-1)^{r(M) - r(A)} (y-1)^{|A| - r(A)}.
\]

By (\ref{eq:chi2}), we have $\chi_M(t) = (-1)^{r(M)} T_M(1-t,0)$.

To put Theorem~\ref{thm:charpolychar} into perspective, we define the {\bf Tutte-Grothendieck ring of matroids} to be the commutative ring $K_0({\rm Mat})$ 
defined as the free abelian group on isomorphism classes of matroids, together with multiplication given by the direct sum of matroids, modulo the relations that if $e$ is neither a loop nor a coloop of $M$ then $[M] = [M \backslash e] + [M / e]$.

If $R$ is a commutative ring, an $R$-valued {\bf Tutte-Grothendieck invariant} is a homomorphism from $K_0({\rm Mat})$ to $R$.
The following result due to Crapo and Brylawski asserts that the Tutte polynomial is the universal Tutte-Grothendieck invariant:

\begin{theorem}
\begin{enumerate}
\item The Tutte polynomial is the unique Tutte-Grothendieck invariant $T : K_0({\rm Mat}) \to \ZZ[x,y]$ satisfying
$T({\rm coloop})=x$ and $T({\rm loop}) = y$.  
\item More generally, if $\phi : K_0({\rm Mat}) \to R$ is any Tutte-Grothendieck invariant then $\phi = \phi_0 \circ T$ where $\phi_0: \ZZ[x,y] \to R$ is the unique ring homomorphism sending
$x$ to $\phi({\rm coloop})$ and $y$ to $\phi({\rm loop})$.
\end{enumerate}
\end{theorem}

Similarly, the characteristic polynomial is the universal Tutte-Grothendieck invariant for combinatorial geometries.  More precisely, if $\phi$ is any Tutte-Grothendieck invariant such that $\phi(M) = \phi(\hat{M})$ for every
loopless matroid $M$, then
\[
\phi(M) = (-1)^{r(M)} \chi_M(1 - \phi({\rm coloop})).
\]

The Tutte polynomial has a number of remarkable properties.  For example, one has the following compatibility with matroid duality:
\[
T_M(x,y) = T_{M^*}(y,x).
\]

\subsection{The rank polynomial} \label{sec:rankpoly}

Let $M$ be a simple matroid with lattice of flats $L$.  The {\bf rank polynomial} of $M$ is
\[
\rho_M(t) = \sum_{F \in L} t^{r(M) - r(F)} = W_0(M)t^{r} + W_1(M) t^{r - 1} + \cdots + W_{r}(M).
\]

The coefficients $W_k(M)$ of $\rho_M(t)$ are strictly positive, and are called the {\bf Whitney numbers of the second kind}.
Concretely, $W_k(M)$ is the number of flats in $M$ of rank $k$.
Comparing with (\ref{eq:chi1}), we see that the coefficients of $\chi_M(t)$ and $\rho_M(t)$ are related by 
\[
\begin{aligned}
w_k(M) &= \sum_{F \in L \; : \; r(F)=k} \mu_L(\emptyset,F), \\
W_k(M) &= \sum_{F \in L \; : \; r(F)=k} 1 .
\end{aligned}
\]

For the matroid $M_n := M(K_n)$ associated to the complete graph $K_n$, $w_k(M_n) = s(n,k)$ and $W_k(M_n) = S(n,k)$ are the Stirling numbers of the first and second kind, respectively (hence the name for the Whitney numbers).

It is conjectured that the Whitney numbers of the second kind form a log-concave, and hence unimodal, sequence for every simple matroid $M$.
This, however, remains an open problem.

It is a recent theorem of Huh--Wang \cite{HuhWang} that if $M$ is a rank $r$ matroid which is representable over some field, then 
$W_1 \leq W_2 \leq \cdots \leq W_{\lfloor r/2 \rfloor}$ and
$W_k \leq W_{r-k}$ for every $k \leq r/2$, see \S\ref{sec:WhitneySecondKind} below.

\subsection{Motivic interpretation of the characteristic polynomial} \label{sec:motivic}
Let $k$ be a field.
The {\bf Grothendieck ring of $k$-varieties} is the commutative ring $K_0({\rm Var}_k)$ 
defined as the free abelian group on isomorphism classes of $k$-varieties, together with multiplication given by the product of varieties, modulo the ``scissors congruence'' relations that whenever $Z \subset X$ is a closed $k$-subvariety we have $[X] = [X \backslash Z] + [Z]$.

When $k=\CC$ or $k=\FF_q$ is a finite field, there is a canonical ring homomorphism\footnote{This homomorphism may be defined when $k=\CC$ as the compactly supported $\chi_y$-genus
from mixed Hodge theory, and when $k=\FF_q$ as 
the compactly supported $\chi_y$-genus in $\ell$-adic cohomology.} $e : K_0({\rm Var}_k) \to \ZZ[t]$ with the property that $e(\AA^1_k)=t$.

Let $A$ be an $r \times m$ matrix with entries in $k$ representing a rank $r$ matroid $M$ with lattice of flats $L$.
Let $V \subset k^m$ be the row space of $A$.
With the notation of Example~\ref{ex:linearflats}, the M\"obius inversion formula shows that in the ring $K_0({\rm Var}_k)$ we have the ``motivic'' identity
\begin{equation} \label{eq:motivic}
[V \cap (k^\times)^m] = \sum_{F \in L} \mu_L(0_L, F) [V \cap L_F].
\end{equation}

(For example, if $V$ is a {\em generic} subspace of $k^m$ then by Inclusion-Exclusion we have
\[
[V \cap (k^\times)^m] = [V \cap L_\emptyset] - \sum_i [V \cap L_i] + \sum_{|I|=2} [V \cap L_I] - \cdots, 
\]
but in general there are subspace relations between the various $V \cap L_I$ governed by the combinatorics of the underlying matroid.)

The identity (\ref{eq:motivic}), which is strongly reminiscent of (\ref{eq:chi1}), can be used to establish Theorem~\ref{thm:charpolychar} in the representable case.  Since $e(V \cap L_F) = t^{r - r(F)}$, it also explains the theorem of Orlik and Solomon \cite{OrlikSolomon} that for $k=\CC$ the Hodge polynomial of 
$V \cap (\CC^*)^m$ is $\chi_M(t)$,
as well as the theorem of Athanasiadis \cite{Athanasiadis} that for $k=\FF_q$ with $q$ sufficiently large we have
\[
|V \cap (\FF_q^\times)^m| = \chi_M(q).
\]

We mentioned in \S\ref{sec:charpoly} that the reduced characteristic polynomial $\bar{\chi}_M(t)$ is the ``projective'' analogue of $\chi_M(t)$.
A concrete way to interpret this statement in the representable case is that since $k^\times$, which satisfies 
$e(k^\times) = t-1$, acts freely on $V \cap (k^\times)^m$, we have
\[
e\left(\PP(V \cap (k^\times)^m)\right) = e(V \cap (k^\times)^m) / e(k^\times) = \chi_M(t) / (t-1) = \bar{\chi}_M(t).
\]

\section{Overview of the proof of the Rota--Welsh Conjecture}

We briefly outline the strategy used by Adiprasito, Huh, and Katz in their proof of the Rota--Welsh conjecture. 
(See \cite{AHKSurvey} for another survey of the proof.)
The first step is to define a {\em Chow ring} $A^*(M)$ associated to an arbitrary loopless matroid $M$. 
The definition of this ring is motivated by work of Feichtner and Yuzvinsky \cite{FY}, who noted that when $M$ is realizable over ${\mathbf C}$, the ring $A^*(M)$ coincides with the usual Chow ring of the de Concini--Procesi ``wonderful compactification'' $Y_M$ of the hyperplane arrangement complement associated to $M$ \cite{DCP,Denham}\footnote{Technically speaking, there are different wonderful compactifictions in the work of de Concini--Procesi; the one relevant for \cite{AHK} corresponds to the ``finest building set''.}.  (Although the definition of $A^*(M)$ is purely combinatorial and does not require any notions from algebraic geometry, it would presumably be rather hard to motivate the following definition without knowing something about the relevant geometric background.)  Note that $Y_M$ is a smooth projective variety of dimension $d := r-1$, where $r$ is the rank of $M$.

\subsection{The Chow ring of a matroid}

Let $M$ be a loopless matroid, and let ${\mathcal F}' = {\mathcal F} \backslash \{ \emptyset,E \}$ be the poset of {\em non-empty proper flats} of $M$. The graded ring $A^*(M)$ is defined as the quotient of the polynomial ring $S_M = {\mathbf Z}[x_F]_{F \in {\mathcal F}'}$ by the following two kinds of relations:
\begin{itemize} 
\item (CH1) For every $a,b \in E$, the sum of the $x_F$ for all $F$ containing $a$ equals the sum of the $x_F$ for all $F$ containing $b$.
\item (CH2) $x_F x_{F'} = 0$ whenever $F$ and $F'$ are incomparable in the poset ${\mathcal F}'$. 
\end{itemize}

The generators $x_F$ are viewed as having degree one. There is an isomorphism\footnote{The isomorphism ${\rm deg}$ should not be confused with the grading on the ring $A^*(M)$, these are two different usages of the term ``degree''.} ${\rm deg} : A^d(M) \rightarrow {\mathbf Z}$ determined uniquely by the property that ${\rm deg}(x_{F_1} x_{F_2} \cdots x_{F_d})=1$ whenever $F_1 \subsetneq F_2 \subsetneq \cdots \subsetneq F_d$ is a maximal flag in ${\mathcal F}'$.

It may be helpful to note that $A^*(M)$ can be naturally identified with equivalence classes of piecewise polynomial functions on the Bergman fan $\Sigma_M$.
The fact that there is a unique homomorphism ${\rm deg} : A^d(M) \rightarrow {\mathbf Z}$ as above means, in the language of tropical geometry, that there is a unique (up to scalar multiple) set of integer weights on the top-dimensional cones of $\Sigma_M$ which make it a {\em balanced} polyhedral complex.

\subsection{Connection to Hodge Theory}

If $M$ is realizable, one can use the so-called {\em Hodge-Riemann relations} from algebraic geometry, applied to the smooth projective algebraic variety $Y_M$ whose Chow ring is $A^*(M)$, to prove the Rota--Welsh log-concavity conjecture for $M$. This is (in retrospect, anyway) the basic idea in the earlier paper of Huh and Katz, about which we will say more in \S\ref{sec:HK} below.

We now quote from the introduction to \cite{AHK}:

\begin{itemize}
\item[ ]  ``While the Chow ring of {M} could be defined for arbitrary {M}, it was unclear how to formulate and prove the Hodge-Riemann relations\ldots We are nearing a difficult chasm, as there is no reason to expect a working Hodge theory beyond the case of realizable matroids. Nevertheless, there was some evidence on the existence of such a theory for arbitrary matroids.''
\end{itemize}

What the authors of \cite{AHK} do is to formulate a purely combinatorial analogue of the Hard Lefschetz Theorem and Hodge-Riemann relations and prove them for the ring $A^*(M)_{\mathbf R} := A^*(M) \otimes {\mathbf R}$ in a purely combinatorial way, making no use of algebraic geometry. The idea is that although the ring $A^*(M)_{\mathbf R}$ is not actually the cohomology ring of a smooth projective variety, from a Hodge-theoretic point of view it behaves as if it were.

\subsection{Ample classes, Hard Lefschetz, and Hodge--Riemann}
In order to formulate precisely the main theorem of \cite{AHK}, we need a combinatorial analogue of hyperplane classes, or more generally of ample divisors. The connection goes through strictly submodular functions.  

A function $c : 2^E \rightarrow {\mathbf R}_{\geq 0}$ is called {\bf strictly submodular}  if  $c(E)=c(\emptyset)=0$ and $c(A \cup B) + c(A \cap B) < c(A) + c(B)$ whenever $A,B$ are incomparable subsets of $E$.  Strictly submodular functions exist, and each such $c$ gives rise to an element $\ell(c) = \sum_{F \in \mathcal{F}'} c(F) x_F \in A^1(M)_\RR.$ The convex cone of all $\ell(c) \in A^1(M)_\RR$ associated to strictly submodular classes is called the {\bf ample cone}\footnote{Actually, the ample cone in \cite{AHK} is 
{\em a priori} larger than what we've just defined, but this subtlety can be ignored for the present purposes.}, and elements of the form $\ell(c)$ are called
{\bf ample classes} in $A^1(M)_\RR$.

Ample classes in $A^1(M)_\RR$ correspond in a natural way to strictly convex piecewise-linear functions on the Bergman fan $\Sigma_M$ (cf.~\S\ref{sec:BergmanFan}).

The main theorem of \cite{AHK} is the following:

\begin{theorem}[Adiprasito--Huh--Katz, 2015] \label{thm:AHKCombHodge}
Let $M$ be a matroid of rank $r=d+1$, let $\ell \in A^1(M)_{\mathbf R}$ be ample, and let $0 \leq k \leq \frac{d}{2}$. Then:
\begin{enumerate}
\item (Poincar{\'e} duality) The natural multiplication map gives a perfect pairing $A^k(M) \times A^{d-k}(M) \rightarrow A^d(M) \cong {\mathbf Z}$.
\item (Hard Lefschetz Theorem) Multiplication by $\ell^{d-2k}$ determines an isomorphism $L_\ell^k : A^k(M)_{\mathbf R} \rightarrow A^{d-k}(M)_{\mathbf R}$.
\item (Hodge-Riemann relations) The natural bilinear form 
\[
Q_\ell^k : A^k(M)_{\mathbf R} \times A^k(M)_{\mathbf R} \rightarrow {\mathbf R}
\] 
defined by $Q_\ell^k(a,b) = (-1)^k a \cdot L_\ell^k b$ is positive definite on the kernel of 
$\ell \cdot L_\ell^k$ (the so-called ``primitive classes'').
\end{enumerate}
\end{theorem}

This is all in very close analogy with analogous results in classical Hodge theory.

\subsection{Combinatorial Hodge theory implies the Rota--Welsh Conjecture}
To see why the Theorem~\ref{thm:AHKCombHodge} implies the Rota--Welsh conjecture, fix $e \in E = \{ 0,\ldots,n \}$. Let $\alpha(e) \in S_M$ be the sum of $x_F$ over all $F$ containing $e$, and let $\beta(e) \in S_M$ be the sum of $x_F$ over all $F$ not containing $e$. 
The images of $\alpha(e)$ and $\beta(e)$ in $A^1(M)$ do not depend on $e$, and are denoted by $\alpha$ and $\beta$, respectively.

\begin{theorem} \label{thm:alphabeta} 
Let $\bar{\chi}_M(t) := \chi_M(t) / (t-1)$ be the reduced characteristic polynomial of $M$, and write $\bar{\chi}_M(t) = m_0 t^d - m_1 t^{d-1} + \cdots + (-1)^d m_d$.  Then $m_k = {\rm deg}(\alpha^{d-k} \cdot \beta^{k})$ for all $k=0,\ldots,d$.
\end{theorem}

The proof of this result is based on the following positive combinatorial formula for $m_k$ due originally to Bj{\"o}rner \cite{BjornerSCM,BjornerMA}.  (It can also be deduced as a straightforward consequence of Weisner's theorem.)

A $k$-step flag $F_1 \subsetneq F_2 \subsetneq \cdots \subsetneq F_k$ in ${\mathcal F}'$ is said to be {\bf initial} if $r_M(F_i)=i$ for all $i$,
and {\bf descending} if 
\[
{\rm min}(F_1) > {\rm min}(F_2) > \cdots > {\rm min}(F_k) > 0,
\]
where for $F \subseteq \{ 0,1,\ldots,n \}$ we set ${\rm min}(F) = {\rm min} \{ i \; : \; i \in F\}$.
 
\begin{proposition} 
$m_k$ is the number of initial, descending $k$-step flags in ${\mathcal F}'$.
\end{proposition}


Although $\alpha$ and $\beta$ are not ample, one may view them as a {\em limit} of ample classes (i.e., they belong to the ``nef cone''). This observation, together with the Hodge-Riemann relations for $A^0(M)$ and $A^1(M)$ and Theorem~\ref{thm:alphabeta}, allows one to deduce the Rota--Welsh conjecture in a formal way.


\subsection{Log-concavity of $f$-vectors of matroids}

The Rota--Welsh conjecture implies a conjecture of Mason and Welsh on $f$-vectors of matroids.

\begin{corollary}[Mason--Welsh Conjecture] \label{cor:MasonWelsh}
Let $M$ be a matroid on $E$, and let $f_k(M)$ be the number of independent subsets of $E$ with cardinality $k$. Then the sequence $f_k(M)$ is log-concave and hence unimodal.
\end{corollary}

To deduce Corollary~\ref{cor:MasonWelsh} from the results of \cite{AHK}, 
one proceeds by showing that the signed $f$-polynomial 
\[
f_0(M) t^r - f_1(M) t^{r-1} + \cdots + (-1)^r f_r(M)
\]
of the rank $r$ matroid $M$ coincides with the reduced characteristic polynomial of an auxiliary rank $r+1$ matroid $M'$ constructed from $M$, the so-called {\em free co-extension of $M$}.\footnote{To define the free co-extension, let $e$ be an auxiliary element not in $E$ and let $E' = E \cup \{ e \}$.  The {\bf free extension of $M$ by $e$} is the matroid $M+e$ on $E'$ whose independent sets are the independent sets of $M$ together with all sets of the form $I \cup \{ e \}$ with $I$ an independent set of $M$ of cardinality at most $r-1$.  The {\bf free co-extension of $M$ by $e$} is the matroid $M \times e$ on $E'$ given by $M \times e = (M^* + e)^*$.}
This identity was originally proved by Brylawski \cite{Brylawski} and subsequently rediscovered by Lenz \cite{Lenz}.

\subsection{High-level overview of the strategy for proving Theorem~\ref{thm:AHKCombHodge}}

The main work in \cite{AHK} is of course establishing Poincar{\'e} duality and especially the Hard Lefschetz Theorem and Hodge--Riemann relations for $M$. 
From a high-level point of view, the proof is reminiscent of Peter McMullen's strategy in \cite{McMullen_Polytopes}, where he 
reduces the so-called ``g-conjecture''\footnote{For an overview of the $g$-conjecture and applications of Hodge theory to the enumerative geometry of polytopes, see e.g. Richard Stanley's article \cite{Stanley}.} for arbitrary simple polytopes to the case of simplices using the ``flip connectivity'' of simple polytopes of given dimension.  

A key observation in \cite{AHK}, motivated in part by McMullen's work, is that for any two matroids $M$ and $M'$ of the same rank on the same ground set $E$, there is a diagram
\[
\xymatrix{\Sigma_{\mathrm{M}} \ar@/^/[r]^{{\rm ``flip"}}& \Sigma_1 \ar@/^/[r]^{{\rm ``flip"}}& \Sigma_2 \ar@/^/[r]^{{\rm ``flip"}} & \cdots \ar@/^/[r]^{{\rm ``flip"}} & \Sigma_{\mathrm{M}'}},
\]


\noindent{w}here each matroidal ``flip''\footnote{A word of caution about the terminology: although these operations are called ÒflipsÓ in \cite{AHK}, they are not analogous to flips in the sense of birational geometry but rather to blowups and blowdowns.} preserves the validity of the Hard Lefschetz Theorem and Hodge-Riemann relations.\footnote{A subtlety is that the intermediate objects $\Sigma_i$ are balanced weighted rational polyhedral fans but not necessarily tropical linear spaces associated to some matroid. So one leaves the world of matroids in the course of the proof, unlike with McMullen's case of polytopes.} 
Using this, one reduces Theorem~\ref{thm:AHKCombHodge} to the Hodge-Riemann relations for projective space, which admit a straightforward (and purely combinatorial) proof.


The inductive approach to the hard Lefschetz theorem and the Hodge-Riemann relations in \cite{AHK} is modeled on the observation that any facet of a permutohedron is the product of two smaller permutohedrons. 


\subsection{Remarks on Chow equivalence} \label{sec:ChowEquiv}
The Chow ring $A^*(M)$ of a rank $d+1$ matroid $M$ on $\{ 0,\ldots,n\}$ coincides with the Chow ring of the smooth but {\em non-complete} toric variety $X(\Sigma_M)$ associated to the Bergman fan of $M$.
One of the subtleties here, and one of the remarkable aspects of the results in \cite{AHK}, is that although
the $n$-dimensional toric variety $X(\Sigma_M)$ is not complete, its Chow ring ``behaves like'' the Chow ring of a $d$-dimensional smooth projective variety.


When $M$ is representable over a field $k$, there is a good reason for this: one can construct a map from a smooth projective variety $Y$ of dimension $d$ to $X(\Sigma_M)$ which induces (via pullback)
an isomorphism of Chow rings
\[
A^*(X(\Sigma_M)) \isomap A^*(Y).
\]
(We call such an isomorphism a {\bf Chow equivalence}.)

For example, if $M = U_{2,3}$ is the uniform matroid represented over $\CC$ by a line $\ell \subset \PP^2$ in general position,
its Bergman fan $\Sigma_M$ is a tropical line in $\RR^2$ (cf.~Example~\ref{ex:tropicalline})
and the corresponding toric variety $X(\Sigma_M)$ is isomorphic to $\PP^2 \backslash \{ 0,1,\infty \}$.
Pullback along the inclusion map $\PP^1 \cong \ell \into \PP^2 \backslash \{ 0,1,\infty \}$ induces a Chow equivalence between 
$\PP^2 \backslash \{ 0,1,\infty \}$ and $\PP^1$.  (However, that the induced map $H^*(\PP^2 \backslash \{ 0,1,\infty \},\CC) \to H^*(\PP^1,\CC)$ on singular cohomology rings is far from being an isomorphism.)

When $M$ is not realizable, however, there is provably no such Chow equivalence between $A^*(M)$ 
and the Chow ring of a smooth projective variety $Y$ mapping to $X(\Sigma_M)$ \cite[Theorem 5.12]{AHK}.

The construction of $Y$ in the realizable case follows from the theory of de Concini--Procesi ``wonderful compactifications''.  
One takes the toric variety $X(\Sigma_{U})$ associated to the $n$-dimensional permutohedron $P_n$ (cf.~\S\ref{sec:BergmanFan}) -- the so-called {\bf permutohedral variety}\footnote{The permutohedral variety is an example of a {\em Losev--Manin moduli space}.} --
and views the Bergman fan $\Sigma_M$ of the realizable rank $d+1$ matroid $M$ as a $d$-dimensional subfan of the normal fan $\Sigma_{U}$ to $P_n$, which is a complete $n$-dimensional fan in $\RR^n$.  This induces an open immersion of toric varieties $X(\Sigma_M) \subset X(\Sigma_{U})$, and the wonderful compactification $Y$ of the hyperplane arrangement complement realizing $M$, which is naturally a closed subvariety of $X(\Sigma_{U})$, belongs to the open subset $X(\Sigma_M)$.  The induced inclusion map $Y \into X(\Sigma_M)$ realizes the desired Chow equivalence.

In this case, 
the linear relations (CH1) come from linear equivalence on the ambient permutohedral toric variety $X(\Sigma_{U})$, pulled back along the open immersion $X(\Sigma_M) \into X(\Sigma_{U})$, and the quadratic relations (CH2) come from 
the fact that if $F$ and $F'$ are incomparable flats then the corresponding divisors are disjoint in $X(\Sigma_U)$.

\subsection{Proof of log-concavity in the realizable case d'apr{\`e}s Huh--Katz}  \label{sec:HK}

The geometric motivation for several parts of the proof of the Rota--Welsh Conjecture comes from the proof of the representable case given in \cite{HK}, and is intimately connected with the geometry of the permutohedral variety.  (We remind the reader, however, that asymptotically 100\% of all matroids are not representable over any field \cite{Nelson}.)
We briefly sketch the argument from \cite{HK}.

The $n$-dimensional permutohedral variety $X(\Sigma_{U})$ is a smooth projective variety which can be considered as an iterated blow-up of $\PP^n$.
After fixing homogenous coordinates on $\PP^n$, we get a number of distinguished linear subspaces of $\PP^n$, for example the $n+1$ points having all but one coordinate equal to zero.
We also get the coordinate lines between any two of those points, and in general we can consider all linear subspaces of the form $\bigcap_{i \in I} H_i$ where $H_i$ is the $i^{\rm th}$ coordinate hyperplane and $I \subset E := \{ 0,1,\ldots,n \}$.
The permtohedral variety $X(\Sigma_{U})$ can be constructed by first blowing up the $n+1$ coordinate points, then blowing up the proper transforms of the coordinate lines, then blowing up the proper transforms of the coordinate planes, and so on.
In particular, this procedure determines a distinguished morphism $\pi_1 : X(\Sigma_{U}) \to \PP^n$ which is a proper modification of $\PP^n$.

There is another distinguished morphism $\pi_2 : X(\Sigma_{U}) \to \PP^n$ which can be obtained by composing $\pi_1$ with the standard Cremona transform ${\rm Crem}: \PP^n \dashrightarrow \PP^n$ given in homogeneous coordinates
by $(x_0:\cdots:x_n) \mapsto (x_0^{-1}:\cdots:x_n^{-1})$.  Although ${\rm Crem}$ is only a rational map on $\PP^n$, it extends to an automorphism of $X(\Sigma_{U})$, i.e., there is a morphism $\widetilde{{\rm Crem}}: X(\Sigma_{U}) \to X(\Sigma_{U})$
such that $\pi_1 \circ \widetilde{{\rm Crem}} = {\rm Crem} \circ \pi_1$ as rational maps $X(\Sigma_{U}) \dashrightarrow \PP^n$.  In other words, $\widetilde{{\rm Crem}}: X(\Sigma_{U}) \to X(\Sigma_{U})$ resolves the indeterminacy locus of ${\rm Crem}$.
We set $\pi_2 = \pi_1 \circ \widetilde{{\rm Crem}}$.

A rank $d+1$ loopless matroid $M$ on $E$ which is representable over $k$ corresponds to a $(d+1)$-dimensional subspace $V$ of $k^{n+1}$ which is not contained in any hyperplane.
Let $\PP(V) \subset \PP^n$ be the projectivization of $V$.  Like $X(\Sigma_U)$ itself, the proper transform $\widetilde{\PP(V)}$ of $\PP(V)$ in $X(\Sigma_{U})$ can be constructed as an iterated blowup, in this case a blowup of $\PP(V)$ at its intersections with the various coordinate spaces of $\PP^n$.  In fact, $\widetilde{\PP(V)}$ coincides with the de Concini--Procesi wonderful compactification $Y$ mentioned above.  The homology class of $\widetilde{\PP(V)}$ in the permutohedral variety depends only on the matroid $M$, and not on the particular choice of the subspace $V$.  We denote by $p_1,p_2$ the restrictions to $\widetilde{\PP(V)}$ of $\pi_1,\pi_2$, respectively.

The key fact from \cite{HK} linking $\widetilde{\PP(V)}$ and the ambient permutohedral variety to the Rota--Welsh Conjecture is the following (compare with Theorem~\ref{thm:alphabeta}):

\begin{theorem} \label{thm:alphabeta2} 
Let $H$ be the class of a hyperplane in ${\rm Pic}(\PP^n)$, let $\alpha = p_1^{-1}(H)$, and let $\beta = p_2^{-1}(H)$.
Then: 
\begin{enumerate}
\item The class of $(p_1 \times p_2)(\widetilde{\PP(V)})$ in the Chow ring of $\PP^n \times \PP^n$ is
\[
m_0 [\PP^d \times \PP^0] + m_1 [\PP^{d-1} \times \PP^1] + \cdots + m_r [\PP^0 \times \PP^d].
\]
\item The $k^{\rm th}$ coefficient $m_k$ of the reduced characteristic polynomial $\bar{\chi}_M(t)$ is equal to ${\rm deg}(\alpha^{d-k} \beta^k)$.
\end{enumerate}
\end{theorem}

The Rota--Welsh conjecture for representable matroids follows immediately from Theorem~\ref{thm:alphabeta2}(2) and the Khovanskii--Teissier inequality, which says that if $X$ is a smooth projective variety of dimension $d$ and $\alpha,\beta$ are nef divisors on $X$ then ${\rm deg}(\alpha^{d-k} \beta^k)$ is a log-concave sequence. 

\subsection{The K{\"a}hler Package}  

The proof of the Khovanskii--Teissier inequality uses Kleiman's criterion to reduce to the case where $\alpha,\beta$ are ample, then uses 
the Kleiman-Bertini theorem to reduce to the case of surfaces, in which case the desired inequality is precisely the classical Hodge Index Theorem.
The Hodge Index Theorem itself is a very special case of the Hodge--Riemann relations.  

One of the original approaches by Huh and Katz to extend their work to non-representable matroids was to try proving a tropical version of the Hodge Index Theorem for surfaces.  However, there are counterexamples to any na{\"i}ve formulation of such a result (see, e.g., \cite[\S{5.6}]{BabaeeHuh}),
and the situation appears quite delicate --- it is unclear what the hypotheses for a tropical Hodge Index Theorem should be and how to reduce the desired inequalities to this special case.  

So instead, inspired by the work of McMullen and Fleming--Karu on Hodge theory for simple polytopes \cite{McMullen_Polytopes,FlemingKaru},
Adiprasito, Huh and Katz developed a completely new method for attacking the general Rota--Welsh Conjecture.

In both the realizable case from \cite{HK} and the general case from \cite{AHK}, one needs only a very special case of the Hodge--Riemann relations to deduce log-concavity of the coefficients of $\bar{\chi}_M(t)$.\footnote{Presumably one can use the general Hodge--Riemann relations to deduce other combinatorial facts of interest about matroids!}  And Poincar{\'e} Duality and the Hard Lefschetz Theorem for Chow rings of matroids are not needed at all for this application.  So it's reasonable to wonder whether Theorem~\ref{thm:AHKCombHodge} is overkill if one just wants a proof of the Rota--Welsh conjecture.  It seems that in practice, Poincar{\'e} Duality, the Hard Lefschetz Theorem, and the Hodge--Riemann relations tend to come bundled together in what is sometimes called the {\em K{\"a}hler package}.\footnote{Note that when some algebraic geometers refer to the K{\"a}hler package, they include additional results such as the Lefschetz hyperplane theorem or K{\"u}nneth formula, which are not part of \cite{AHK}.}  This is the case, for example, in the 
algebro-geometric work of de Cataldo--Migliorini and Cattani
\cite{DCM,Cattani},
in the work of McMullen and Fleming--Karu on Hodge theory for simple polytopes \cite{McMullen_Polytopes,FlemingKaru}, in the work of Elias--Williamson \cite{EliasWilliamson} on Hodge theory for Soergel bimodules, and in Adiprasito--Huh--Katz. 

In the case of simple polytopes and the $g$-conjecture, what is needed is in fact the Hard Lefschetz Theorem, and not the Hodge--Riemann relations, for the appropriate Chow ring.  But again the proof proceeds by establishing the full K{\"a}hler package.

One of the important differences between \cite{AHK} and \cite{FlemingKaru}, already mentioned above, is that the intermediate objects in the inductive procedure from \cite{AHK}, obtained by applying flips to Bergman fans of matroids, are no longer themselves Bergman fans of matroids (whereas in \cite{FlemingKaru} all of the simplicial fans which appear come from simple polytopes).  
Another important difference is that in
the polytope case one is working with $n$-dimensional fans in $\RR^n$,
whereas in the matroid case one is working with $d$-dimensional fans in $\RR^n$, where $d<n$ except in the trivial (but important) case of the $n$-dimensional permutohedral fan.  In both the polytope and matroid situations the fan in question defines an $n$-dimensional toric variety, but the toric variety is projective in the polytope case and non-complete in the matroid case.
As mentioned above in \S\ref{sec:ChowEquiv}, the ``miracle'' in the matroid case is that the Chow ring of the $n$-dimensional non-complete toric variety $X(\Sigma_M)$ behaves as if it were the Chow ring of a $d$-dimensional smooth projective variety; in particular, it satisfies Poincar{\'e} Duality, Hard Lefschetz, and Hodge--Riemann of ``formal'' dimension $d$.

\subsection{Whitney numbers of the second kind} \label{sec:WhitneySecondKind}

The Whitney numbers of the second kind $W_k(M)$ (cf.~\S\ref{sec:rankpoly}) are much less tractable then their first-kind counterparts.  
In particular, the log-concavity conjecture for them remains wide open.  However, there has been recent progress by Huh and Wang \cite{HuhWang} concerning a related conjecture, the so-called ``top-heavy conjecture'' of Dowling and Wilson:

\begin{conjecture} \label{conj:DowlingWilson}
Let $M$ be a matroid of rank $r$.  Then for all $k < r/2$ we have $W_k(M) \leq W_{r-k}(M)$.
\end{conjecture}

In analogy with the work of Huh--Katz, Huh and Wang prove:

\begin{theorem}[Huh--Wang, 2016] \label{thm:HuhWang}
For all matroids $M$ representable over some field $k$:
\begin{enumerate}
\item The first half of the sequence of Whitney numbers of the second kind is unimodal, i.e., $W_1(M) \leq W_2(M) \leq \cdots \leq W_{\lfloor r/2 \rfloor}(M)$.
\item Conjecture~\ref{conj:DowlingWilson} is true.
\end{enumerate}
\end{theorem}

The following corollary is a generalization of the de Bruijn-Erd\H{o}s theorem that every non-collinear set of points $E$ in a projective plane determines at least $|E|$ lines:

\begin{corollary}
Let $V$ be a $d$-dimensional vector space over a field and let $E$ be a subset which spans $V$.  Then (in the partially ordered set of subspaces spanned by subsets of $E$), there are at least as many $(d-k)$-dimensional subspaces as there are $k$-dimensional subspaces, for every $k \leq d/2$.  \end{corollary}

We will content ourselves with just a couple of general remarks concerning the proof of Theorem~\ref{thm:HuhWang}.
Unlike in the Rota--Welsh situation of Whitney numbers of the first kind, the projective algebraic variety $Y_M'$ which one associates to $M$ in this case is 
highly singular; thus instead of invoking the K{\"a}hler package for smooth projective varieties, Huh and Wang have to use analogous but much harder results about intersection cohomology.
Specifically, they require the Bernstein--Beilinson--Deligne--Gabber decomposition theorem for intersection complexes\footnote{See \cite{CDM_Bulletin} for an overview of the decomposition theorem and its many applications.} and the Hard Lefschetz theorem for $\ell$-adic intersection cohomology of projective varieties.

It is tempting to fantasize about a proof of Conjecture~\ref{conj:DowlingWilson} along the lines of \cite{AHK}.
One of many significant challenges in this direction would be to 
construct a combinatorial model for intersection cohomology of the variety $Y_M'$.





\bibliographystyle{amsplain}
\bibliography{CEB}


\end{document}